\documentclass[a4paper,12pt]{article}

\usepackage{amsmath,amssymb,amsthm}

\usepackage{mathrsfs}

\makeatother
\usepackage{graphicx}
\newtheorem{theorem}{Theorem}
\newtheorem{definition}{Definition}
\newtheorem{remark}{Remark}

\begin{document}

\title{Optimal Control Problem in a Stochastic Model\\
with Periodic Hits on the Boundary of a Given Subset\\
of the State Set (Tuning Problem)}

\author{V.~P.~Shnurkov\thanks{National Research University
``Higher School of Economics",
Moscow Institute of Electronics and Mathematics,
Department of Applied Mathematics,
pshnurkov@hse.ru}}

\date{}

\maketitle

\begin{abstract}
In this paper, a general stochastic model with controls applied at the moments
when the random process hits the boundary of a given subset of the state set
is proposed and studied.
The general concept of the model is formulated
and its possible applications in technical and economic systems are described.
Two versions of the general stochastic model,
the version based on the use of a continuous-time semi-Markov process
with embedded absorbing Markov chain
and the version based on the use of a discrete-time Markov process
with absorbing states,
are analyzed.
New representations of the stationary cost index of the control quality
are obtained for both versions.
It is shown that this index can be represented
as a linear-fractional integral functional of two discrete probability distributions
determining the control strategy.
The results obtained by the author of this paper about an extremum of such functionals
were used to prove that, in both versions of the model,
the optimal control is deterministic and is determined
by the extremum points of functions of two discrete arguments
for which the explicit analytic representations are obtained.
The perspectives of the further development of this stochastic model
are outlined.
\smallskip

\end{abstract}

\section*{Introduction}

The following phenomenon is observed when many processes of stochastic character
encountered in technical and economic systems are analyzed.
At certain times, the basic stochastic process describing the system under study
can leave a given subset of states which are assumed to be admissible.
The process can be returned to the subset of admissible states
by an external action which can be treated as a control in the stochastic model.
As a result of the control action, the process returns to one of the states
in the admissible set.
In this case, the state itself (or the number of this state)
into which the process is taken by the external action is a control parameter or a decision.
If the decision is accepted by the stochastic scheme,
then it is assumed that the state into which the process is taken
is determined by a certain probability distribution which
can be called a controlling distribution.
After the basic stochastic process is transferred to one of the admissible states,
it again begins to evolve independently of the past, i.e.,
its behavior depends only on the state where it was taken by the
external action (control).
At a certain time, the basic process again leaves the set of admissible states,
and then the external control action is applied according to the above-cited rules.
The problem of control optimization is to find controlling probability distributions
on which a certain parameter of the control efficiency
attains an extremum.

The above-formulated ideas, which are based on specific characteristics of operation
of actual systems,
form the concept of a stochastic model with controls applied at the moments
when the basic stochastic process attains the boundary of a given subset of the state set.
Let us consider the possible applications of such a model
in technical and economic systems in more detail.

For a technical system, a subset of admissible states
can be interpreted as the set of operation states,
and the output from this set means the system failure.
The transition of the basic process into the subset of admissible states
(external control action)
is a reconstructing action applied to this system
to re-establish its operation.
The state of the basic process immediately after the control action
determines the form of the reconstructing action, i.e.,
the character of the reconstruction.
Thus, in applications,
the problem of choosing the optimal control action
in the form of optimal controlling probability distribution
is the problem of determining the optimal form of reconstruction
after the system failure.

The corresponding stochastic control model
is also meaningful for economic systems.
The basic stochastic process is assumed to describe
the variations in some economic parameter
(price of financial assets on the stock market,
prices of a product on the goods market,
currency price at currency auctions).
The set of admissible states can be interpreted
as the set of values of this parameter,
which are favorable or acceptable for a special agent of the market
applying an external control action
to the corresponding financial parameter.
Such a special agent of the market
is usually an authorized governmental structure.
The output of the basic process from a given subset of admissible states
is undesirable from the standpoint of this market agent.
The control action can involve different actions of the market agent
aimed at changing the values of the corresponding financial parameter.
A specific example of such actions
is the central bank currency intervention
aimed at changing the exchange rate of a currency.
A similar example is the grain intervention performed by the grain fund,
an authorized structure of the Ministry of Agriculture,
aimed at changing the prices of different cereal crops.
As a results of such actions,
the process value (basic parameter) returns to the set of admissible states
and the evolution of this parameter value goes on
independently of the past
and depends only on the new initial value
formed by the external control action.
The above-cited financial parameter (basic process)
continues to evolve according to the laws of free market
in the same economical situation as at the preceding state
until it again leaves the subset of admissible states
and a new control action is applied.

It follows from the above considerations that
developing a stochastic model with controls at the moments of
hitting the boundary of a certain subset of the state set
and solving the corresponding optimal control problem
is an actual problem of applied mathematics.
In what follows, such a problem is called a \textbf{tuning problem},
which reflects the above-described external action
on the basic parameter of a technical or an economic system.

In the present paper, we consider the following two possible versions
of the general stochastic model discussed above.
The first version (continuous time)
is based on a semi-Markov stochastic process
with finitely many states, which has an embedded absorbing Markov chain.
The absorbing states are assumed to the boundary states,
and the other states are internal and admissible.
After absorption, a control action
taking the process into one of the internal admissible states
is applied.
The further evolution of the process continues independently of the past
and is described by the probability characteristics
of the initial semi-Markov process until the next absorption occurs.
The second version (discrete time) involves
the Markov chain with two absorbing states
which are assumed to be boundary states.
The other states are assumed to be internal and admissible.
The evolution of the discrete-time model
obeys similar laws.
In both versions,
we pose and solve optimal control problems
for stationary cost indices which,
according to their economic content,
are the mean specific profits.
Thus, two tuning problems related to various stochastic models
are solved in this paper.

\section{Extremum problem for linear-fractional\\
integral functionals defined on a set\\
of discrete probability distributions}

To solve the optimal control problem formulated below,
we first recall a certain general result of the theory of extremum problems,
namely, a theorem about an unconditional extremum
for a linear-fractional integral functional on a set of probability measures.
Such a theorem was formulated by the author in~[1],
where it was also shown that this result is a theoretical foundation
for solving optimal control problems
for a semi-Markov stochastic process with a finite set of states.
Further,
the theorem about an unconditional extremum
for a linear-fractional integral functional was proved in~[2]
and its application to solving the optimal control problem for a semi-Markov process
was justified there in detail.
We note that the stochastic scheme of control
at the moments of hitting the boundary of a given subset of the state set,
which is proposed in this paper,
is a special version of general control of a semi-Markov process,
which differs from the version considered in fundamental works [3,4]
and in contemporary studies
(see, e.g., [5,6]).
Thus, the tuning problem does not reduce to the standard control problem
for a semi-Markov process,
and a separate analysis is required to solve this problem.
To explain such an analysis, we consider a special version of assertions
about an unconditional extremum of a linear-fractional integral functional
defined on a set of discrete probability distributions.

We consider a set of discrete sets $U_i={1,2,\dots,n_i}$, $i=1,2,\dots,N$, $N<\infty$.
The number of elements in each set~$U_i$ can be either finite or countable:
$n_i\leq\infty$, $i=1,2,\dots,N$.
In what follows, these sets are interpreted as sets of admissible solutions (controls)
accepted in different states of the stochastic model,
but in this section, they are abstract.
On each set $U_i$, we introduce a collection
of all possible probability distributions of the form
$\alpha^{(i)}=(\alpha_1^{(i)},\alpha_2^{(i)},\dots,\alpha_{n_i}^{(i)})$,
$\alpha_s^{(i)}\geq 0$, $s=1,2,\dots,n_i$;
$\sum\limits_{s=1}^{n_i}\alpha_s^{(i)}=1$, $i=1,2,\dots,N$.
We denote such a set of probability distributions defined on~$U_i$
by~$\Gamma_{d}^{(i)}$, $i=1,2,\dots,N$.
Further, we consider the Cartesian product of spaces
$U=U_1\times U_2\times\dots\times U_N$.
Following the classical scheme of introducing the measure
on a Cartesian product of spaces~[7],
we introduce the probability measure on~$U$
as the product of probability measures on the spaces
$U_1,U_2,\dots,U_N$ determined by the distributions
$\alpha^{(1)},\alpha^{(2)},\dots,\alpha^{(N)}$.
Thus, the probability measure on~$U$ is given
by the set of probability distributions
$\alpha^{(1)}$, $\alpha^{(2)},\dots,\alpha^{(N)}$.
We denote the set of probability measures on~$U$ by $\Gamma_{d}$.
Several additional conditions related to the set~$\Gamma_{d}$
are described below in the statement of the extremum problem.

Similarly to [1,2], we introduce the notion of degenerate discrete probability distribution.

\begin{definition}\rm %1
A probability distribution $\alpha^{(i)*}(k_i)$
is said to be \textit{degenerate} if
$\alpha_{k_i}^{(i)}=1$, $\alpha_l^{(i)}=0$, $l=1,2,\dots,n_i$, $l\neq k_i$.
A point $k_i\in U_i$ is called a \textit{point of concentration}
of a degenerate distribution $\alpha^{(i)*}(k_i)$.
As is known, a degenerate distribution corresponds to
the deterministic quantity taking the value~$k_i$.
\end{definition}

We denote the set of degenerate probability distributions defined on~$U_i$ by
$\Gamma_{d}^{(i)*}$, $i=1,2,\dots,N$.
Obviously, there is a one-to-one correspondence between the sets
$\Gamma_d^{(i)*}$ and $U_i$, $i=1,2,\dots,N$.
We accordingly denote
the set of all degenerate probability measures defined on the set~$U$
by~$\Gamma_{d}^{*}$.
Each degenerate measure in the set $\Gamma_{d}^{*}$
is defined by a set of degenerate distributions
$\alpha^{(1)*},\alpha^{(2)*},\dots,\alpha^{(N)*}$.

We assume that two numerical functions are defined on the set~$U$:
$$
A(k_1,k_2,\dots,k_N): U\rightarrow R,\quad
B(k_1,k_2,\dots,k_N): U\rightarrow R,\quad
\text{where $k_i\in U_i$, $i=1,2,\dots,N$}.
$$

We note that the integral over a discrete measure defined on a discrete set
can be transformed into a sum,
and the corresponding multidimensional integral over the measure
generated by the product of initial measures
defined on the Cartesian product of discrete spaces
becomes a multidimensional sum.
According to this, we introduce the following definition by analogy with~[1,2].

\begin{definition}\rm %2
A \textit{linear-fractional integral functional} (\textit{in the discrete version})
or simply a \textit{discrete linear-fractional integral functional}
defined on a set of collections of discrete probability distributions $\Gamma_{d}$
is defined to be the mapping
$I(\alpha^{(1)},\alpha^{(2)},\dots,\alpha^{(N)}): \Gamma_{d}\rightarrow R$
given by the expression
\begin{equation}
I(\alpha^{(1)},\alpha^{(2)},\dots,\alpha^{(N)})=\dfrac
{\sum\limits_{k_1=1}^{n_1}\sum\limits_{k_2=1}^{n_2}\dots\sum\limits_{k_{N}=1}^{n_N}
A(k_1,k_2,\dots,k_{N})\alpha_{k_1}^{(1)}\alpha_{k_2}^{(2)}\dots\alpha_{k_{N}}^{(N)}}
{\sum\limits_{k_1=1}^{n_1}\sum\limits_{k_2=1}^{n_2}\dots\sum\limits_{k_{N}=1}^{n_N}
B(k_1,k_2,\dots,k_{N})\alpha_{k_1}^{(1)}\alpha_{k_2}^{(2)}\dots\alpha_{k_{N}}^{(N)}}
\end{equation}
\end{definition}

\begin{definition}\rm %3
A function $C(k_1,k_2,\dots,k_N):U\rightarrow R$ defined by the expression
\begin{equation}
C(k_1,k_2,\dots,k_N)=\dfrac{A(k_1,k_2,\dots,k_N)}{B(k_1,k_2,\dots,k_N)},
\end{equation}
is the \textit{test} function of the discrete linear-fractional integral functional
$I(\alpha^{(1)},\alpha^{(2)},\dots,\alpha^{(N)})$ given by formula~(1).
\end{definition}

Let us formulate the corresponding extremum problem for
$I(\alpha^{1},\alpha^{(2)},\dots,\alpha^{(N)})$ of the form (1)
on a set of collections of discrete probability distributions $\Gamma_d$:
\begin{equation}
I(\alpha^{1},\alpha^{(2)},\dots,\alpha^{(N)}) \to \mathop{extr}, \qquad
(\alpha^{1},\alpha^{(2)},\dots,\alpha^{(N)})\in \Gamma_d.
\end{equation}

We assume that some preliminary conditions similar to the corresponding conditions
introduced when solving the extremum problem
for a linear-fractional integral functional of general structure considered in~[1,2]
are satisfied for the above-posed extremum problem~(3).
Let us specify these conditions.
\begin{enumerate}
\item
The functionals of discrete probability distributions,
which determine the numerator and denominator in expression~(1),
are defined for any probability distributions
$(\alpha^{(1)},\alpha^{(2)},\dots,\alpha^{(N)})\in \Gamma_d$ as
\begin{equation}
I_1(\alpha^{(1)},\alpha^{(2)},\dots,\alpha^{(N)})=
\sum\limits_{k_1=1}^{n_1}\sum\limits_{k_2=1}^{n_2}\dots\sum\limits_{k_{N}=1}^{n_N}
A(k_1,k_2,\dots,k_{N})\alpha_{k_1}^{(1)}\alpha_{k_2}^{(2)}\dots\alpha_{k_{N}}^{(N)},
\end{equation}
\begin{equation}
I_2(\alpha^{(1)},\alpha^{(2)},\dots,\alpha^{(N)})=
\sum\limits_{k_1=1}^{n_1}\sum\limits_{k_2=1}^{n_2}\dots\sum\limits_{k_{N}=1}^{n_N}
B(k_1,k_2,\dots,k_{N})\alpha_{k_1}^{(1)}\alpha_{k_2}^{(2)}\dots\alpha_{k_{N}}^{(N)}.
\end{equation}
In other words, the numerical series in the right-hand sides of expressions~(4) and~(5)
are assumed to converge.
\item
For any discrete probability distributions
$(\alpha^{(1)},\alpha^{(2)},\dots,\alpha^{(N)})\in\Gamma_d$,
the functional $I_2(\alpha^{(1)},\alpha^{(2)},\dots,\alpha^{(N)})$
does not vanish, namely,
\begin{equation*}
I_2(\alpha^{(1)},\alpha^{(2)},\dots,\alpha^{(N)})=
\sum\limits_{k_1=1}^{n_1}\sum\limits_{k_2=1}^{n_2}\dots\sum\limits_{k_{N}=1}^{n_N}
B(k_1,k_2,\dots,k_{N})\alpha_{k_1}^{(1)}\alpha_{k_2}^{(2)}\dots\alpha_{k_{N}}^{(N)}\ne 0.
\end{equation*}
\item
The set of collections of degenerate probability distributions $\Gamma_d^*$
is completely contained in the set $\Gamma_d:\Gamma_d^*\subset\Gamma_d$.
\end{enumerate}

\begin{remark}\rm %1
As in the general version~[1],  conditions~2 and~3 imply that
the function $B(k_1,k_2,\dots,k_N)$ is strictly of constant sign
for all $(k_1,k_2,\dots,k_N)\in U$.
At the same time,
if this condition related to the character of the function $B(k_1,k_2,\dots,k_N)$
is satisfied, then condition~2 is satisfied automatically.
In~[2], it is specially noted that the condition
of being strictly of constant sign
(and, in particular, of being strictly positive)
for the function $B(k_1,k_2,\dots,k_N)$ is natural for
optimal control problems for semi-Markov processes.
In this connection, it is required that this condition is satisfied
in the fundamental theorem on the solution of extremum problem~(3).
\end{remark}

\begin{definition}\rm %4
The set of collections of discrete probability distributions $\Gamma_{d}$
is said to be \textit{admissible} in extremum problem~(3)
if conditions~1 and~3 in the system of preliminary conditions
are satisfied.
\end{definition}

We now formulate the fundamental theorem on the solution of extremum problem~(3)
which is a particular case of Theorem~1 formulated in~[1].
We restrict our consideration to the first assertion in this theorem
as the most important for solving the optimal control problem
considered in the present paper.

\begin{theorem} %1
Assume that the set of collections of discrete probability distributions
$\Gamma_{d}$ in extremum problem~{\rm(3)} is admissible
and the function $B(k_1,k_2,\dots,k_N)$
in discrete linear-fractional integral functional~{\rm(1)}
is strictly of constant sign
{\rm(}strictly positive or strictly negative{\rm)}
for all values of the arguments $(k_1,k_2,\dots,k_N)\in U$.
Assume also that the test function
of the discrete linear-fractional integral functional
$C(k_1,k_2,\dots,k_{N})=\dfrac{A(k_1,k_2,\dots,k_{N})}{B(k_1,k_2,\dots,k_{N})}$
attains a global extremum {\rm(}maximum or minimum{\rm)}
on the set~$U$ at a fixed point $(k_1^*,k_2^*,\dots,k_{N}^*)$.
Then the solution of the corresponding extremum problem~{\rm(3)}
for the maximum or minimum exists
and is attained on the set of degenerate probability distributions
$(\alpha^{(1)^*},\alpha^{(2)^*},\dots,\alpha^{(N)^*})$
concentrated at the respective points $k_1^*,k_2^*,\dots,k_{N}^*$
and the following relations are satisfied:
\begin{align}
\max\limits_{(\alpha^{(1)},\alpha^{(2)},\dots,\alpha^{(N)})\in\Gamma_d}
&I(\alpha^{(1)},\alpha^{(2)},\dots,\alpha^{(N)})
=
\max\limits_{(\alpha^{(1)^*},\alpha^{(2)^*},\dots,\alpha^{(N)^*})\in\Gamma_d^*}
I(\alpha^{(1)^*},\alpha^{(2)^*},\dots,\alpha^{(N)^*})
\nonumber\\
&=\max\limits_{(k_1,k_2,\dots,k_{N})\in U}
\dfrac{A(k_1,k_2,\dots,k_{N})}{B(k_1,k_2,\dots,k_{N})}
= \dfrac{A(k_1^*,k_2^*,\dots,k_{N}^*)}{B(k_1^*,k_2^*,\dots,k_{N}^*)}
\end{align}
if the global maximum of the test function is attained
at the point $(k_1^*,k_2^*,\dots,k_{N}^*)$;
\begin{align}
\min\limits_{(\alpha^{(1)},\alpha^{(2)},\dots,\alpha^{(N)})\in\Gamma_d}
&I(\alpha^{(1)},\alpha^{(2)},\dots,\alpha^{(N)})
= \min\limits_{(\alpha^{(1)^*},\alpha^{(2)^*},\dots,\alpha^{(N)^*})\in\Gamma_d^*}
I(\alpha^{(1)^*},\alpha^{(2)^*},\dots,\alpha^{(N)^*})
\nonumber\\
&=\min\limits_{(k_1,k_2,\dots,k_{N})\in U}
\dfrac{A(k_1,k_2,\dots,k_{N})}{B(k_1,k_2,\dots,k_{N})}=
\dfrac{A(k_1^*,k_2^*,\dots,k_{N}^*)}{B(k_1^*,k_2^*,\dots,k_{N}^*)}
\end{align}
if the global minimum of the test function is attained
at the point $(k_1^*,k_2^*,\dots,k_{N}^*)$.
\end{theorem}

\section{General structure of semi-Markov model\\
with controls applied at the moments\\
of hitting the boundary}

We first assume that all stochastic objects introduced below
are defined on a certain initial probability space
$(\Omega, \mathscr{A}, P)$
formalizing an actual random experiment.
The general structure of the probability space
is in detail described in fundamental works ([8], Vol.~1; [9].

The proposed continuous-time model is based on
the semi-Markov stochastic process with a discrete state set.
The general theory of such processes is presented
in the classical work~[10] and in contemporary studies
(see, e.g.,~[11]).
We also note that a concise description of the foundations of the theory
of semi-Markov processes can be found in the book~[9].

Let $\xi^{(n)}(t)$, $n=0,1,2,\dots$, be a sequence of independent semi-Markov processes
with absorption and equal probability characteristics.
By $X=\{0,1,\dots,N\}$, $N<\infty$, we denote the set of states of these processes.

We additionally agree about the character of the set of states.
It is natural to assume that~$\{0\}$ and~$\{N\}$
are boundary states for the state set~$X$,
and $\{1,2,\dots,N-1\}$ are internal states.
We assume that the set $\{1,2,\dots,N-1\}$ is admissible
in the framework of our control model,
and the inputs into boundary states are outputs from an admissible set.
We assume that the boundary states are absorbing, i.e.,
the times of input into these states are moments of absorption.
To make the application of the theory of processes with absorptions
more convenient, we redenote the states in the initial set
$\{0,1,\dots,N\}$ as follows:
the state $\{0\}$ is still denoted by $\{0\}$,
the state $\{N\}$ is now denoted by $\{1\}$,
and the other states $\{1,2,\dots,N-1\}$ are denoted by $\{2,3,\dots,N\}$,
respectively.
Thus, in the state set $X=\{0,1,\dots,N\}$,
the states~$\{0\}$ and~$\{1\}$ are boundary and absorbing,
and the states $\{2,3,\dots,N\}$ are internal and admissible.

In what follows, we assume that
necessary initial probability characteristics are prescribed
for the semi-Markov processes $\xi^{(n)}(t)$, $n=0,1,2,\dots$.
Different forms of such characteristics are given in~[10].
In particular, one can introduce
semi-Markov functions $Q_{ij}(t)$, $i,j\in X$,
which are joint distributions
of the embedded Markov chain transitions
and durations of the process stay in different states.

We also introduce two systems of independent nonnegative random variables
$\{\Delta_n^{(0)}$, $n=0,1,2,\dots\}$ and $\{\Delta_n^{(1)}$, $n=0,1,2,\dots\}$
whose distributions for each fixed~$n$
depend on an additional condition related to this model.

We assume that in the structure of the proposed model,
the semi-Markov processes $\xi^{(n)}(t)$, $n=0,1,2,\dots$,
describe the evolution of the initial system
on the periods between the control actions.
Each process $\xi^{(n)}(t)$ starts to evolve
in one of the admissible states $i\in\{2,3,\dots,N\}$.
In a finite time after the evolution beginning,
the process $\xi^{(n)}(t)$ is in one
of the absorbing boundary states $\{0\}$ or $\{1\}$
with probability~$1$.
After absorption of the process $\xi^{(n)}(t)$,
the model experiences the so-called control action,
which is manifested in the transition from the boundary state
into one of the admissible (internal) states.
The probabilities of the transition from the boundary state $0$
into internal admissible states are given by the vector
$\alpha^{(0)}=\left(\alpha_k^{(0)},~k=2,3,\dots,N\right)$,
$\sum\limits_{k=2}^N\alpha_k^{(0)}=1$,
and the probabilities of the transition from the boundary state~$1$
into internal admissible states are given by the vector
$\alpha^{(1)}=\left(\alpha_k^{(1)},~k=2,3,\dots,N\right)$,
$\sum\limits_{k=2}^N\alpha_k^{(1)}=1$.
The probability distributions~$\alpha^{(0)}$ and~$\alpha^{(1)}$
describe the external control actions.
The time of this transition is a random variable $\Delta_n^{(0)}$
if the transition occurs from the state~$\{0\}$
and a random variable $\Delta_n^{(1)}$
if the transition occurs from the state~$\{1\}$.
In this case,
the distribution of the random variable $\Delta_n^{(0)}$ or $\Delta_n^{(1)}$
can generally depend on the number of the admissible state $k \in \{2,3,\dots,N\}$
into which the transition leads.
After the transition into an admissible state $k \in \{2,3,\dots,N\}$,
the evolution of the system is independent of the past
(for a fixed initial state~$k$)
and is described by the semi-Markov process $\xi^{(n+1)}(t)$
whose probability characteristics coincide with
the corresponding characteristics of the process $\xi^{(n)}(t)$.
The further evolution is similar.
%The possible trajectory of the above-described continuous-time semi-Markov model
%is shown in Fig.~1.

\begin{figure}[h]
\centerline{
    \includegraphics[width=0.7\textwidth]{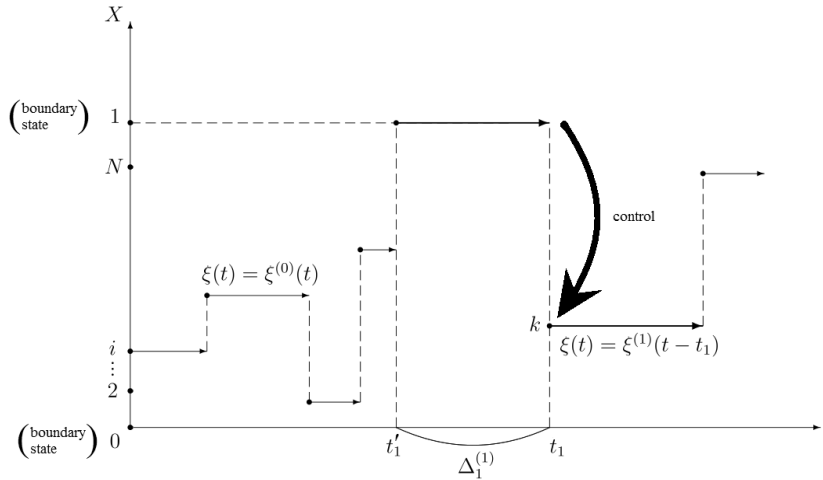}
}
\caption{Trajectory of the semi-Markov process $\xi(t)$
controlled at the moments of hitting the boundary of a given subset of the state set.}
\end{figure}

The corresponding discrete-time stochastic model is constructed similarly.
Instead of the semi-Markov process $\xi^{(n)}(t)$, $n=0,1,2,\dots$,
we use independent absorbing Markov chains
$\xi^{(n)}=\left\{\xi_t^{(n)},~t=0,1,2,\dots\right\}$, $n=0,1,2,\dots$,
ranging in the above-defined state set $X$
with the same transition probability matrix $\textbf{P}$.
We note that, in the discrete-time model,
the durations of a single stay in all states
can conditionally be assumed to be~$1$.
Therefore, the random variables $\{\Delta_n^{(0)}$, $n=0,1,2,\dots\}$
and $\{\Delta_n^{(1)}$, $n=0,1,2,\dots\}$
describing the control action duration
are also formally equal to~$1$.

The optimal control problem in both versions of the above-described model
is formally posed as the extremum problem on the set of pairs of
discrete probability distributions
$\left(\alpha^{(0)},\alpha^{(1)}\right)$
with respect to a certain stationary cost index of efficiency
which has the meaning of the mean specific profit.

\section{General probability characteristics of the model}

In this section, we determine all probability and cost characteristics
required to solve the optimal control problems
in both versions of the model under study.
We divide these characteristics into initial or primary,
which must be given to describe the model,
and secondary,
which can be expressed in terms of the primary ones
using the well-known theoretical assertions.
We begin with determining the initial (primary) characteristics.

We let $\textbf{P}$ denote the transition probability matrix
for the absorbing Markov chain embedded in the semi-Markov process
$\xi^{(n)}(t)$, $n=0,1,2,\dots$.
As is known in the theory of absorbing Markov chains~[12],
such a matrix has the following block structure
$$
\textbf{P}=\left(\;\begin{matrix} \textbf{P}_{11} & \textbf{P}_{10} \cr
\textbf{P}_{01} & \textbf{P}_{00} \cr \end{matrix}\;\right),
$$
where $\textbf{P}_{00}$
is the transition probability matrix inside the set
of irreversible admissible states $\{2,3,\dots,N\}$;

$\textbf{P}_{01}$ is the probability matrix for the transitions
from the set of admissible states $\{2,3,\dots,N\}$
into the set of absorbing boundary states $\{0,1\}$;

$\textbf{P}_{10}$ is the zero matrix corresponding
to transitions from absorbing states $\{0,1\}$
into the set of admissible states $\{2,3,\dots,N\}$;

$\textbf{P}_{11}$ is the unit matrix corresponding
to transitions inside the set of absorbing states $\{0,1\}$.

For the discrete-time model, we preserve the above-introduced notation
for the transition probability matrices
for the independent absorbing Markov chains $\xi^{(n)}$, $n=0,1,2,\dots$.

For the other primary characteristics, we use two types of notation
distinguishing the corresponding characteristics of the discrete-time model
by the symbol ``$~\widehat{~}~$'' above a letter.
Let us introduce these characteristics in succession:

$\tau_k$ is the mathematical expectation (mean) of the time of stay
of the semi-Markov process $\xi^{(n)}(t)$ in the state $k\in\{2,3,\dots,N\}$.
According to the above remarks,
there are no such quantities in the discrete-time model;

$c_k$ is the mathematical expectation of the income obtained
during the total time of a single stay
of the semi-Markov process $\xi^{(n)}(t)$ in the state $k\in\{2,3,\dots,N\}$;
in the discrete-time model,
the corresponding quantities are denoted by $\widehat{c}_k$, $k\in\{2,3,\dots,N\}$,
and they characterize the value of the income obtained during a single stay
of the Markov chain in the state $k\in\{2,3,\dots,N\}$;

$d_k^{(0)}$ and $d_k^{(1)}$ are the mathematical expectations
of costs due to the control action,
i.e., to the transfer of the basic process from the boundary states $\{0,1\}$
into a state $k\in\{2,3,\dots,N\}$.
In the discrete-time model, the corresponding characteristics are
denoted by $\widehat{d}_k^{~(0)}$ and $\widehat{d}_k^{~(1)}$,
$k\in \{ 2,3,\dots,N \}$;
according to the economic content of these quantities,
we assume that they are negative;

$\mu_k^{(0)}$ and $\mu_k^{(1)}$ are mathematical expectations
of the time necessary to apply the control action,
i.e., to transfer the basic process from the boundary states $\{0,1\}$
into a state $k\in\{2,3,\dots,N\}$.
In the discrete-time model, such characteristics are absent,
because such transitions are realized at one step of the model.

Now let us determine the secondary characteristics.
To express these characteristics in terms of primary ones,
we use the theoretical results for Markov chains with absorption
mainly following the work~[12].

By $\textbf{M}=\left(m_{ij}\right)$ we denote the matrix given by the formula
$\textbf{M}=\left(\textbf{I}-\textbf{P}_{00}\right)^{-1}$.
The entries of this matrix $m_{ij}$ are
the mathematical expectations of the number of inputs
of the Markov chain embedded in the semi-Markov process $\xi^{(n)}(t)$
into a state $j\in\{2,3,\dots,N\}$
during the time till the moment of absorption,
which are determined under the condition that this chain starts to evolve
from a state $i\in\{2,3,\dots,N\}$.
The corresponding mathematical expectations for the discrete-time model
are denoted by $\widehat{m}_{ij}$, $i,j\in\{2,3,\dots,N\}$.
These quantities are entries
of the matrix $\widehat{\textbf{M}}=\left(\widehat{m}_{ij}\right)$
given by a similar formula.

Let $m_i$ be the mathematical expectation of the time
from the beginning of the process $\xi^{(n)}(t)$ evolution
to the moment of absorption,
which is determined under the condition that the initial state
is~$i$.
The corresponding mathematical expectation in the discrete-time model
is denoted by~$\widehat{m}_i$, $i\in\{2,3,\dots,N\}$.
These quantities are given by the formulas
\begin{equation}
m_i=\sum\limits_{j=2}^N m_{ij}\tau_j,~i\in\{2,3,\dots,N\},
\end{equation}
\begin{equation}
\widehat{m}_i=\sum\limits_{j=2}^N \widehat{m}_{ij},~i\in\{2,3,\dots,N\}.
\end{equation}

By $r_i$ we denote the mathematical expectation of the income
due to the evolution of the process $\xi^{(n)}(t)$
till the absorption,
which is determined under the condition that the initial state
is~$i$, $i\in\{2,3,\dots,N\}$.
The corresponding mathematical expectation in the discrete-time model
is denoted by~$\widehat{r}_i$, $i\in\{2,3,\dots,N\}$.
These quantities are determined by the formulas~([4,12])
\begin{equation}
r_i=\sum\limits_{j=2}^N m_{ij}c_j,~i\in\{2,3,\dots,N\};
\end{equation}
\begin{equation}
\widehat{r}_i=\sum\limits_{j=2}^N \widehat{m}_{ij}\widehat{c}_j,~i\in\{2,3,\dots,N\}.
\end{equation}

Further, we let $b_{i0}$ and $b_{i1}$ denote the conditional probabilities
of the events that the semi-Markov process $\xi^{(n)}(t)$
is absorbed in the respective state~$0$ or~$1$,
which are determined under the condition that the initial state of the process
is $i\in\{2,3,\dots,N\}$.
We let the vectors
$(b_{i0}$ and $b_{i1})$, $i\in\{2,3,\dots,N\}$,
form a $(N-1)\times 2$ matrix $\textbf{B}$.

Then we have the formula~[12]
\begin{equation}
\textbf{B}=(\textbf{I}-\textbf{P}_{00})^{-1}\textbf{P}_{01}=\textbf{MP}_{01}.
\end{equation}

The corresponding probabilities for the discrete-time model are denoted by
$\widehat{b}_{i0}$ and $\widehat{b}_{i1}$, $i\in\{2,3,\dots,N\}$.
They are entries of the matrix
$\widehat{\textbf{B}}=\left(\widehat{b}_{i0},~\widehat{b}_{i1},~~i\in\{2,3,\dots,N\}\right)$,
which is given by a similar formula
\begin{equation}
\widehat{\textbf{B}}=\left(I-P_{00}\right)^{-1}P_{01}=\widehat{\textbf{M}}P_{01}.
\end{equation}

Thus, the secondary probability characteristics of continuous- and discrete-time models,
given by formulas (8)--(13),
can explicitly be expresses in terms of the initial (primary) characteristics.
This means that to solve the optimal control problems in these models,
it suffices to prescribe only a set of primary characteristics.

\section{Analytic representation of stationary\\
cost indices of the control quality}

Now we consider the stationary cost index of the control efficiency
related to the evolution of this stochastic model.
We note that the nature of such an index is the same in both versions of the model,
with continuous time and with discrete time.
From the formal standpoint, this index is the limit
$$
I=\lim\limits_{t\rightarrow\infty}\dfrac{V(t)}{t},
$$
where $V(t)$ is the mathematical expectation of an additive functional
(accumulation functional) defined on the trajectory of the stochastic process under study.

According to its economic content, the index $I$
is the mean specific profit obtained as
the system operates under stationary conditions.
Such indices have been known since the fundamental studies
in the theory of control of Markov and semi-Markov stochastic processes~[3,4].
We note that similar indices are used in contemporary studies~[5,6].

We use ergodic theorems for Markov and semi-Markov stochastic processes
to obtain explicit analytic representations for such indices in both versions
of the stochastic model considered in this paper.
For continuous-time models, one can use the assertion proved in~[4], Chapter~5, Section~5.5,
or the corresponding assertions for general semi-Markov models with arbitrary state spaces
and controls~[5,6].
For discrete-time models,
one can use the classical ergodic theorem for Markov chains ([13], Chapter~8]).

\begin{remark}\rm %2
To apply the above-cited ergodic theorems,
it is necessary to require that
the constructed continuous- and discrete-time stochastic models
have a certain property which can be called stability.
For continuous-time models,
this means that the probabilities of absorption for any initial state
$i\in\{2,3,\dots,N\}$ must satisfy the conditions
$b_{i0}>0$, $b_{i1}>0$, and $b_{i0}+b_{i1}=1$.
Then the absorbing boundary states are admissible
from any initial state $i$ with probability~$1$.

Now we consider an auxiliary Markov chain with the state set $\{0,1\}$,
which is formed by the values of the process at the moments of hitting the boundary states
$\xi'_n=\xi(t'_n)$, $n=0,1,2,\dots$.
Under these conditions, this chain is irreducible
and has a unique stationary distribution.
[The classical theory of Markov chains is discussed in detail in [8], Chapter~VIII;
also see the above-cited work [13], Chapter~8.]
At the same time, the additive cost functional $V(t)$
determining the control efficiency index~$I$
can be defined on the trajectories of precisely this Markov chain.
According to the above-cited ergodic theorems,
the limits $I=\lim\limits_{t\rightarrow\infty}\big[\frac{V(t)}{t}\big]$
exist.
All quantities entering the explicit representation for~$I$
are expressed in terms of the initial (primary and secondary) characteristics of the model.
Similar remarks also remain true for the discrete-time model.
\end{remark}

Now we formulate assertions about analytic representations of stationary cost indices
of the control efficiency in terms of the well-known primary and secondary probability characteristics
of the model and the discrete probability distributions $\alpha^{(0)}$ and $\alpha^{(1)}$
describing the control actions.
The following assertion holds for the continuous-time model.

\begin{theorem} %2
The stationary index of the mean specific profit
has the analytic representation
\begin{equation}
I=
\dfrac
{\sum\limits_{k=2}^N\alpha_k^{(0)}\left[d_k^{(0)}+r_k\right]\sum\limits_{l=2}^N\alpha_l^{(1)}b_{l0}+
\sum\limits_{k=2}^N\alpha_k^{(1)}\left[d_k^{(1)}+r_k\right]\sum\limits_{l=2}^N\alpha_l^{(0)}(1-b_{l0})}
{\sum\limits_{k=2}^N\alpha_k^{(0)}\left[\mu_k^{(0)}+m_k\right]\sum\limits_{l=2}^N\alpha_l^{(1)}b_{l0}+
\sum\limits_{k=2}^N\alpha_k^{(1)}\left[\mu_k^{(1)}+m_k\right]\sum\limits_{l=2}^N\alpha_l^{(0)}(1-b_{l0})}.
\end{equation}
\end{theorem}

The corresponding assertion for the discrete-time model can be formulated as follows.

\begin{theorem} % 3.
The stationary index of the mean specific profit can be represented as
\begin{equation}
\widehat{I}=\dfrac{\sum\limits_{k=2}^N \alpha_k^{(0)}\left[\widehat{d}_k^{~(0)}+\widehat{r}_k\right]
\sum\limits_{l=2}^N \alpha_l^{(1)}\widehat{b}_{l0}
+\sum\limits_{k=2}^N \alpha_k^{(1)} \left[\widehat{d}_k^{~(1)}+\widehat{r}_k\right]
\sum\limits_{l=2}^N \alpha_l^{(0)}(1-\widehat{b}_{l0})}
{\sum\limits_{l=2}^N \alpha_l^{(1)}\widehat{b}_{l0}+\sum\limits_{l=2}^N \alpha_l^{(0)}(1-\widehat{b}_{l0})}.
\end{equation}
\end{theorem}

\begin{remark} \rm %3
The continuous- and discrete-time stochastic models under study
can be generalized by assuming that the set of boundary states
is an arbitrary finite set.
In this case, the states contained in it are called external states
with respect to the internal admissible states.
The auxiliary Markov chain formed by the values of the process
at the moments of transition to external states
has finitely many states.
Under natural assumptions similar to those formulated in Remark~2,
this chain is irreducible and has a unique stationary distribution.
The control actions are determined by finitely many probability distributions
each of which described the transitions from a certain external state
into an internal (admissible) state.
For such models,
one can prove assertions similar to Theorems~2 and~3.
In this case, expressions~(14) and~(15) contain stationary probabilities
of the auxiliary Markov chains mentioned above.
In turn, these stationary probabilities are solutions of systems of linear algebraic equations
and can be expressed in terms of the initial probability characteristics of stochastic models.
\end{remark}

\section{Solution of the optimal control problem}

Now to complete the study, we present the solution of the optimal control problems
in both versions of the stochastic model.
For this, we represent the stationary cost indices of the control quality,
given by formulas~(14) and~(15),
as linear-fractional integral functionals of discrete probability distributions
$\alpha^{(0)}$ and $\alpha^{(1)}$.
We use analytic transformations of formula (14)
to prove the following assertion for the continuous-time model.

\begin{theorem} %4
The stationary cost index $I(\alpha^{(0)},\alpha^{(1)})$ given by formula~{\rm(14)}
can be represented as a discrete linear-fractional integral functional
of the probability distributions $\alpha^{(0)},\alpha^{(1)}$, namely,
\begin{equation}
I(\alpha^{(0)},\alpha^{(1)})=\dfrac{\sum\limits_{l_0=2}^N
\sum\limits_{l_1=2}^NA(l_0,l_1)\alpha_{l_0}^{(0)}\alpha_{l_1}^{(1)}}
{\sum\limits_{l_0=2}^N\sum\limits_{l_1=2}^NB(l_0,l_1)\alpha_{l_0}^{(0)}\alpha_{l_1}^{(1)}},
\end{equation}
\begin{equation}
\text{where}~A(l_0,l_1)=\left[d_{l_0}^{(0)}+r_{l_0}\right]b_{l_1,0}+\left[d_{l_1}^{(1)}+r_{l_1}\right]b_{l_0,1},
\end{equation}
\begin{equation}
B(k_0,k_1)=\left[\mu_{l_0}^{(0)}+m_{l_0}\right]b_{l_1,0}+\left[\mu_{l_1}^{(1)}+m_{l_1}\right]b_{l_0,1}.
\end{equation}
\end{theorem}

It follows from Theorem~4 that the test function of the discrete linear-fractional integral functional~(16)
is given by the relation
\begin{equation}
C(l_0,l_1)=\dfrac{A(l_0,l_1)}{B(l_0,l_1)},
\end{equation}
where the functions $A(l_0,l_1)$ and $B(l_0,l_1)$
are respectively determined by formulas~(17) and~(18).

A similar assertion, which can be proved by using analytic transformations of formula~(15),
also holds for the discrete-time model.

\begin{theorem} % 5
The stationary cost index
$\widehat{I}\left(\alpha^{(0)},\alpha^{(1)}\right)$ given by formula~{\rm(15)}
can be represented as a discrete linear-fractional integral functional
of the probability distributions $\alpha^{(0)}$ and $\alpha^{(1)}$, namely,
\begin{equation}
\widehat{I}\left(\alpha^{(0)},\alpha^{(1)}\right)
=\dfrac{\sum\limits_{l_0=2}^N\sum\limits_{l_1=2}^N \widehat{A}(l_0,l_1)\alpha_{l_0}^{(0)}\alpha_{l_1}^{(1)}}
{\sum\limits_{l_0=2}^N\sum\limits_{l_1=2}^N \widehat{B}(l_0,l_1)\alpha_{l_0}^{(0)}\alpha_{l_1}^{(1)}},
\end{equation}
where
\begin{equation}
\widehat{A}(l_0,l_1)=\left[\widehat{d}_{l_0}^{~(0)}+\widehat{r}_{l_0}\right]\widehat{b}_{l_1,0}
+\left[\widehat{d}_{l_1}^{~(1)}+\widehat{r}_{l_1}\right]\widehat{b}_{l_0,1},
\end{equation}
\begin{equation}
\widehat{B}(l_0,l_1)=\widehat{b}_{l_1,0}+\widehat{b}_{l_0,1}.
\end{equation}
\end{theorem}

It follows from the assertion of Theorem~5 that the test function
of the discrete linear-fractional integral functional~(20)
has the form
\begin{equation}
\widehat{C}(l_0,l_1)=\dfrac{\widehat{A}(l_0,l_1)}{\widehat{B}(l_0,l_1)},
\end{equation}
where the functions $\widehat{A}(l_0,l_1)$ and $\widehat{B}(l_0,l_1)$
are determined by formulas~(21) and~(22), respectively.

The optimal control problems for continuous- and discrete-time models
are formulated in the form of extremum problems~(3)
posed on the set $\Gamma_{d}$
which consists of pairs of all possible discrete probability distributions
$\left(\alpha^{(0)},~\alpha^{(1)}\right)$
with respect to functionals~(16) and~(20), respectively.
In this case, the sets of admissible controls~$U_0$ and~$U_1$
are determined as $U_0=U_1=\{2,3,\dots,N\}$,
and the set~$U$ is the Cartesian product
$U=U_0\times U_1 = \left\{(l_0,l_1),~l_0\in U_0,~l_1\in U_1\right\}$.
To solve these problems, it is necessary to use Theorem~1.
We note that the conditions of Theorem~1 are satisfied in both problems.
In particular,
the functions $B(l_0,l_1)$ and $\widehat{B}(l_0,l_1)$
are strictly positive for all $(l_0,l_1)\in U$.

If the set of admissible controls $U_0=U_1=\{2,3,\dots,N\}$ is finite,
then the test functions of the discrete linear-fractional integral functionals
$C(l_0,l_1)$ and $\widehat{C}(l_0,l_1)$, given by formulas~(19) and~(23),
attain the global extrema (maximum and minimum).
Then, it follows from Theorem~1 that the solution of the optimal control problems
exists for any of the model.
The optimal controls for both models are deterministic
and are determined by the points at which the test functions
attain their global maxima.
If the set of admissible controls $U_0=U_1=\{2,3,\dots,N\}$ is countable ($N=\infty$),
then the solutions of extremum problems for continuous- and discrete-time models
are also completely determined by the extremum properties of the test functions
$C(l_0,l_1)$ and $\widehat{C}(l_0,l_1)$,
as it follows from the general assertion about an extremum
of the linear-fractional integral functional ([1], Theorem~1).
Thus, this theorem is a theoretical foundation which permits obtaining
the complete solution of the tuning problem for the stochastic models considered above.

\section*{Concluding remarks}

In the study whose main results are described in the present paper,
we developed a general concept of stochastic model with controls
applied at the moments when the basic process hits
the boundary of a given subset of the state set.
The proposed concept is realized for two versions
of the initial stochastic model:
the version based on the use of a continuous-time semi-Markov process
with embedded absorbing Markov chain
and the version based on the use of a discrete-time Markov process
with absorbing states.
For both versions of the model, we pose and solve
the optimal control problems, i.e., the corresponding tuning problems.
The theoretical background of such solutions is the theorem
about an unconditional extremum of a linear-fractional integral functional~[1],
which, in the discrete version, becomes Theorem~1 in Section~1.

We note that the possible constructions of new versions of the general stochastic model
with controls at the moments of hitting the boundary of a given subset of the state set
are far from being complete.
We distinguish several fundamental conditions necessary to construct
such well-posed and meaningful versions.
\begin{enumerate}
\item
The set of states must contain a given subset of the so-called admissible states
and a finite subset of states which we conditionally call external or boundary states
with respect to the admissible states.
In this case, the set of admissible states need not be discrete.
\item
A certain general condition related to the behavior of the initial stochastic process
$\xi^{(n)}(t)$ (continuous time) or $\xi^{(n)}$ (discrete time)
must be satisfied.
This condition means that this stochastic process,
which starts to evolve from an arbitrary admissible state,
must get into one of the external (boundary) states in a finite time
with probability~$1$.
\item
The Markov property must be satisfied
at the moments of transition into external (boundary) states
and at the moments of output from external states
and transition into one of the states in the admissible set.
\item
The probability characteristics describing the transition
from an arbitrary fixed admissible state
into one of the external states (an analog of the absorption probabilities)
can be determined analytically.
\item
The mathematical expectations of the time passed
from leaving an arbitrary fixed admissible state and entering one of the external states
and the mathematical expectations of the income accumulated in this time
can be determined analytically.
\end{enumerate}

We note that conditions~1--3 in this set of conditions
are related to the general stochastic properties of the model.
Conditions~4 and~5 must ensure the possibility of determining
the analytic representations for the required probability characteristics of the basic process
which are related to the time of its stay in the set of admissible states.

Under the above-listed conditions,
one can construct a stochastic model with controls
similar to the models which were considered in this study
and obtain an explicit representation for the stationary cost index of the control efficiency
and the complete solution of the corresponding tuning problem.


\begin{thebibliography}{99}

\bibitem{sh1}
V.~P.~Shnurkov,
``Solution of the unconditional extremum problem
for a linear-fractional integral functional on a set of probability measures,''
Dokl. Ross. Akad. Nauk \textbf{470} (4), 387--392 (2016)
[Dokl. Math.].

\bibitem{sh2}
V.~P.~Shnurkov, A.~K.~Gorshenin, and V.~V.~Belousov,
``Analytic solution of optimal control problem for semi-Markov process with a finite set of states,''
Informatika i Primeneniya \textbf{10} (4), 72--88 (2016).

\bibitem{sh3}
W.~S.~Jewell,
``Markov-Renewal Programming. I, II.''
Operations Research \textbf{11} (6), 938--971 (1963).

\bibitem{sh4}
H.~Mine and S.~Osaki,
\textit{Markovian Decision Processes}
(Amer. Elsevier, New York, 1970).

\bibitem{zh5}
F.~Luque-Vasquez and O.~Herndndez-Lerma,
``Semi-Markov control models with average costs,''
Appl. Math. \textbf{26} (3), 315--331 (1999).

\bibitem{sh6}
O.~Vega-Amaya and F.~Luque-Vasquez,
``Sample-path average cost optimality for semi-Markov control processes on Borel spaces:
Unbounded costs and mean holding times,''
Appl. Math. \textbf{27} (3), 343--367 (2000).

\bibitem{sh7}
P.~Halmos,
\textit{Measure Theory}
(Litton, New York, 1950; Springer-Verlag, New York, 1974).

\bibitem{sh8}
A.~N.~Shiryaev,
\textit{Probability}
(Springer, 2016).

\bibitem{sh9}
A.~A.~Borovkov,
\textit{Probability Theory}
(Springer. 2013).

\bibitem{sh10}
V.~S.~Korolyuk and A.~F.~Turbin,
\textit{Semi-Markov Processes and Their Applications}
(Naukova Dumka, Kiev, 1976)
[in Russian].

\bibitem{sh11}
J.~Janssen and R.~Manca,
\textit{Applied Semi-Markov Processes}
(Springer, New York, 2006).

\bibitem{sh12}
J.~Kemeny and J.~Snell,
\textit{Finite Markov Chains}
(Prentice--Hall, 1959).

\bibitem{sh13}
V.~S.~Korolyuk, N.~I.~Portenko, A.,~V.~Skorokhod, and A.~F.~Turbin,
\textit{Handbook in Probability and Mathematical Statistics},
Ed. by V.~S.~Korolyuk
(Naukova Dumka, Kiev, 1995)
[in Russian].

\end{thebibliography}
\end{document}